\documentclass[12pt]{article}
\usepackage{amsmath}
\usepackage{amsfonts}
\usepackage{amssymb}
\newtheorem{theorem}{Theorem}
\newtheorem{definition}[theorem]{Definition}

\begin{document}

\title{\textbf{Plane partitions and their pedestal polynomials\thanks{%
Part of this work has been carried out in the framework of the Labex
Archimede (ANR-11-LABX-0033) and of the A*MIDEX project
(ANR-11-IDEX-0001-02), funded by the \textquotedblleft Investissements
d'Avenir" French Government programme managed by the French National
Research Agency (ANR). Part of this work has been carried out at IITP RAS.
The support of Russian Foundation for Sciences (project No. 14-50-00150) is
gratefully acknowledged.}}}
\author{Oleg Ogievetsky$^{\ddag }$\thanks{%
On leave of absence from P. N. Lebedev Physical Institute, Leninsky Pr. 53,
117924 Moscow, Russia} $\,$ and Senya Shlosman$^{^{\ddag },^{\ddag \ddag }}$
\\
$^{^{\ddag }}$Aix Marseille Universit\'{e}, Universit\'{e} de Toulon, \\
CNRS, CPT UMR 7332, 13288, Marseille, France\\
$^{^{\ddag \ddag }}$Inst. of the Information Transmission Problems,\\
RAS, Moscow, Russia }
\date{}
\maketitle

\begin{abstract}\noindent
We define, for an arbitrary partially ordered set, a multi-variable polynomial generalizing the hook polynomial.

\end{abstract}

\section{Introduction}

Let ${\cal{S}}$ be a partially ordered set. In this work we associate to ${\cal{S}}$ a multi-variable polynomial $\mathfrak{h}$.
When ${\cal{S}}$ is a Young diagram, the principal specialization of $\mathfrak{h}$ coincides with the hook polynomial.

Our construction of $\mathfrak{h}$ begins with defining a polynomial $\mathfrak{h}_P$, where $P$ is an arbitrary linear extension of ${\cal{S}}$.
Then we show that in fact $\mathfrak{h}_P$ does not depend on $P$. The proof
uses the equality  (\ref{01}) (precise definitions are given in sections \ref{Mainresult} and \ref{Proofanddiscussion}), which is implied by the bijection
between the set of reverse partitions on ${\cal{S}}$ and the product of the set of $P$-pedestals on ${\cal{S}}$ 
and the set of Young diagrams with at most $\vert {\cal{S}}\vert$ rows.
It would be interesting to find a direct, not referring to the formula (\ref{01}), proof of the theorem \ref{thmnotdep}.

To facilitate the exposition we take for ${\cal{S}}$ the set of nodes of a Young diagram $\lambda$. In this situation, linear extensions of ${\cal{S}}$
correspond to standard Young tableaux of shape $\lambda$, see Definition \ref{getopo}. Our results and proofs work
in the same way for general ${\cal{S}}$ (and linear extensions of ${\cal{S}}$ instead of standard Young tableaux).

\section{Main result}\label{Mainresult}

Let $\lambda=(\lambda_{1},\dots,\lambda_{l})\vdash n$, $\lambda_{1}\geq \dots\geq\lambda_{l}>0$, be a partition of $n$, $\lambda_{1}+\dots+\lambda
_{l}=n$. We identify $\lambda$ with its Young diagram, that is, the set of nodes
\[\alpha=(i,j)\ \ \text{with}\ \ j=1,\dots,\lambda_{i}\ \ \text{for each}\ \ i=1,\dots,l.\]

A \textit{standard Young tableau of shape} $\lambda$ is a bijection $Q\colon\lambda\rightarrow\{1,\dots,n\}$ such that the function $Q(i,j)$
increases in $i$ and $j$. We denote the set of these standard Young tableaux by $\mathfrak{st}_{\lambda}.$

\begin{definition}\label{getopo}
\label{noteidenflin}\hspace{-0.3cm} .\hspace{0.2cm}\textrm{Let $\preccurlyeq$ be
the minimal transitive partial order on $\lambda$, containing $(i,j)\prec (i+1,j)$ and $(i,j)\prec (i,j+1)$ for all
possible $(i,j)$. A standard Young tableau of shape $\lambda$ can be
identified with a linear extension of $\preccurlyeq$, that is, a linear order, compatible with the partial order $\preccurlyeq$.}
We denote by $\preccurlyeq_P$ the linear order associated to a standard tableau $P$.
\end{definition}

Let $\mathbb{Z}_{\geqslant0}$ be the set of non-negative integers. A \textit{reverse plane partition of shape} $\lambda$ is a function
$\mathfrak{Q}\colon\lambda\rightarrow\mathbb{Z}_{\geqslant0}$, non-decreasing in $i$ and $j$. It is \textit{column-strict} if it increases in $j$. We
visualize reverse plane partitions by placing the number $\mathfrak{Q}(i,j)$ in the node $(i,j)$ for each $(i,j)\in\lambda$.
We denote by $\left\vert
\mathfrak{Q}\right\vert $ the volume of $\mathfrak{Q}$,
$\left\vert \mathfrak{Q}\right\vert
=\sum_{(i,j)\in\lambda}\mathfrak{Q}(i,j).$

Let $\lambda$ be a partition of $n$. Let ${\mathcal{\bar{S}}}_{\lambda}$ be the set of reverse plane partitions
of shape $\lambda,$ and ${\mathcal{S}}_{\lambda}\subset{\mathcal{\bar{S}}}_{\lambda}$ the subset of reverse column-strict plane partitions
of shape $\lambda$. Recall that the Schur function $s_{\lambda}$ is the formal power series in infinitely many variables $\mathbf{x}=\left(
x_{0},x_{1},x_{2},...\right)  ,$ given by
\[s_{\lambda}\left(  \mathbf{x}\right)  =\sum_{\mathfrak{Q\in}{\mathcal{S}
}_{\lambda}}\prod_{\alpha\in\lambda}x_{\mathfrak{Q}\left(  \alpha\right)  }.\]
We need the similarly defined formal power series $\bar{s}_{\lambda}$,
\[\bar{s}_{\lambda}\left(  \mathbf{x}\right)  =\sum_{\mathfrak{Q\in
}{\mathcal{\bar{S}}}_{\lambda}}\prod_{\alpha\in\lambda}x_{\mathfrak{Q}\left(
\alpha\right)  }.\]
The Schur function $s_{\lambda}\left(  \mathbf{x}\right)$, unlike our `wrong' Schur function, $\bar{s}_{\lambda}\left(  \mathbf{x}\right)  ,$
is symmetric (see, e.g. \cite{St}).

To define the pedestal polynomial we proceed as follows. Let $P,Q\in \mathfrak{st}_{\lambda}$ be two standard Young tableaux of shape $\lambda$.
We are going to compare the corresponding linear orders $\preccurlyeq_P$ and $\preccurlyeq_Q$ on the set of the nodes of $\lambda$.
Let $\alpha_{1}=\left(  1,1\right)  \prec_{Q}\alpha_{2}\prec_{Q}\dots\prec_{Q}\alpha_{n}$ be the list of all the nodes of $\lambda,$ enumerated
according to the order $\preccurlyeq_Q$. We say that a node $\alpha_{k}$ is a $\left(  P,Q\right)  $-disagreement node, if $\alpha_{k+1}\prec_{P}\alpha_{k}$
(while, of course, $\alpha_{k}\prec_{Q}\alpha_{k+1}$). We define the reverse plane partition $q_{_{P,Q}}$ of shape $\lambda$ by
\begin{equation}q_{_{P,Q}}\left(  \alpha_{k}\right)  =\sharp\left\{  l\ :\ l<k\ ,\ \alpha
_{l}\text{ is a }\left(  P,Q\right)  \text{-disagreement node}\right\}.\label{13}\end{equation}
Indeed, the function $q_{_{P,Q}}$ is non-decreasing with respect to the order $\preccurlyeq_Q$, hence $q_{_{P,Q}}$ is a reverse plane partition.

The reverse plane partition $q_{_{P,Q}}$ thus defined is called $P$-pedestal of $Q$, see \cite{S}. We finally define the polynomial
$\mathfrak{h}_{P}\left(\mathbf{x}\right)  $ by
\begin{equation}
\mathfrak{h}_{P}\left(  \mathbf{x}\right)  =\sum_{Q\in\mathfrak{st}_{\lambda}
}\prod_{\alpha\in\lambda}x_{q_{_{P,Q}}\left(  \alpha\right)  }. \label{02}
\end{equation}
The set (as $Q$ runs over $\mathfrak{st}_{\lambda}$) of $P$-pedestals depends on $P$. However, the polynomial $\mathfrak{h}_{P}$ has the following remarkable property. 

\begin{theorem}
\label{thmnotdep} \hspace{-.3cm} .\hspace{.2cm} The function $\mathfrak{h}
_{P}\left(  \mathbf{x}\right)  $ does not depend on $P$ from $\mathfrak{st}
_{\lambda}.$
\end{theorem}

\begin{definition}
\hspace{-.3cm} .\hspace{.2cm} We call the function
\[\mathfrak{h}_{\lambda}\left(  \mathbf{x}\right)  =\mathfrak{h}_{P}\left(
\mathbf{x}\right)  ,\]
where $P$ is any standard Young tableau $P$ of shape $\lambda$, the \textbf{pedestal polynomial}.
\end{definition}

Our theorem states that the pedestal polynomial is well-defined. For example, $\mathfrak{h}_{(3,2)}\left(  \mathbf{x}\right)=x_0^5+x_0^4x_1+x_0^3x_1^2+x_0^2x_1^3+
x_0^2x_1^2x_2$.

\section{Proof and discussion}\label{Proofanddiscussion}

\vskip.2cm \noindent\textit{Proof of the Theorem.} Let $R_{n}$ be the component of degree $n$ of the ring of formal power series in variables
$\mathbf{x}=\left(  x_{0},x_{1},x_{2},...\right)  $ with integer coefficients. The element $u\left(  \mathbf{x}\right)  $ of $R_{n}$ is a sum,
\[u\left(  \mathbf{x}\right)  =\sum_{i_{1}\leq i_{2}\leq...\leq i_{n}}
a_{i_{1}i_{2}...i_{n}}\mathbf{x}_{i_{1},i_{2},\dots,i_{n}}\ ,\ \ \ \mathbf{x}
_{i_{1},i_{2},\dots,i_{n}}:=x_{i_{1}}x_{i_{2}}...x_{i_{n}},\]
where $i_{1},i_{2},...,i_{n}$ are non-negative integers, and the coefficients $a_{i_{1}i_{2}...i_{n}}$ are integer.
For example, the functions $s_{\lambda}\left(\mathbf{x}\right)$ and $\bar{s}_{\lambda}\left(  \mathbf{x}\right)  $ belong to $R_{n}$. The
function $\bar{s}_{(n)}$, corresponding to the one-row partition $\lambda=(n)$,
\[\bar{s}_{(n)} \left(  \mathbf{x}\right)  =\sum_{i_{1}\leq i_{2}\leq...\leq i_{n}
}\mathbf{x}_{i_{1},i_{2},\dots,i_{n}}\in R_{n}\]
will play a special role in our argument.

We define the $\ast$-product on monomials by
\[\mathbf{x}_{i_{1},i_{2},\dots,i_{n}}\ast\mathbf{x}_{j_{1},j_{2},\dots,j_{n}
}:=\mathbf{x}_{i_{1}+j_{1},i_{2}+j_{2},\dots,i_{n}+j_{n}}\]
and extend it by linearity to the ring structure on $R_{n}$. The ring $(R_{n},\ast)$ is isomorphic to a subring
of the ring $\mathbb{C}\left[\left[ y_{1},...,y_{n}\right]\right]$ of formal power series in $n$ variables, via
the monomorphism $\varphi_{n}:R_{n}\to \mathbb{C}\left[  \left[ y_{1},...,y_{n}\right]  \right]$,
defined by
\[\varphi_{n}\left(  \mathbf{x}_{i_{1},i_{2},\dots,i_{n}}\right)  =y_{1}^{i_{1}}...y_{n}^{i_{n}}.\]
In particular, $(R_{n},\ast)$ inherits from $\mathbb{C}\left[\left[y_{1},...,y_{n}\right]\right]$ the property of having no zero divisors.

Fix a standard Young tableau $P\in\mathfrak{st}_{\lambda}$. We will prove now the identity
\begin{equation}\bar{s}_{\lambda}\left(  \mathbf{x}\right)  =\mathfrak{h}_{P}\left(\mathbf{x}\right)  \ast \bar{s}_{(n)} \left(  \mathbf{x}\right)  , \label{01}\end{equation}
which implies, due to the absence of zero divisors in $(R_{n},\ast)$, the assertion of the theorem, since the first and the last terms in the identity do not depend on $P$.

The bijective proof of the identity $\left(  \ref{01}\right)  $ follows from \cite{S}. Relations (46), (48) and
$(50)$ of that paper describe bijections $b^{\phantom{-1}}_{St}$, $b_{St}^{-1}$ between the set of reverse plane partitions of shape $\lambda$
and the product of the set of $P$-pedestals and the set of Young diagrams with at most $n$ rows.
Let $\mathfrak{Q}$ be a reverse plane partition and $b^{\phantom{-1}}_{St}(\mathfrak{Q})=(q,\mu)$, where $q$ is a $P$-pedestal and $\mu$ a Young diagram.
The construction of $b^{\phantom{-1}}_{St}$ (see below) implies that the monomial, corresponding to $\mathfrak{Q}$ in
$\bar{s}_{\lambda}\left(\mathbf{x}\right)$ is the $\ast$-product of the monomial corresponding to $q$ in $\mathfrak{h}_{P}\left(\mathbf{x}\right)$ and the monomial corresponding to $\mu$
in $\bar{s}_{(n)} \left(  \mathbf{x}\right)$, and the proof of (\ref{01})  follows. These bijections were in fact used earlier by D. Knuth in \cite{K}.

The bijections $b^{\phantom{-1}}_{St}$, $b_{St}^{-1}$ are defined as follows.

Let $\mathfrak{Q}$ be a reverse plane partition of shape $\lambda.$ Then we can define the partition of $\left\vert \mathfrak{Q}\right\vert $ with at most
$n$ rows,
\begin{equation}p=\Pi\left(  \mathfrak{Q}\right)  \label{u12},\end{equation}
by just listing all the entries of the two-dimensional array of values of $\mathfrak{Q}$ in the non-increasing order.

To every reverse plane partition $\mathfrak{Q}$ we associate the standard Young tableau $Q\left(  \mathfrak{Q}\right)  \in\mathfrak{st}_{\lambda}$ as
follows. Define the linear order $\preccurlyeq_{\mathfrak{Q}}$ on the nodes of $\lambda$ by
\[\alpha^{\prime}\prec_{\mathfrak{Q}}\alpha^{\prime\prime}\text{ if
}\mathfrak{Q}\left(  \alpha^{\prime}\right)  <\mathfrak{Q}\left(
\alpha^{\prime\prime}\right)  \text{ or if }\mathfrak{Q}\left(  \alpha
^{\prime}\right)  =\mathfrak{Q}\left(  \alpha^{\prime\prime}\right)  \text{
and }\alpha^{\prime}\prec_{P}\alpha^{\prime\prime}.\]
Then $Q\left(  \mathfrak{Q}\right)  $ is defined by the relation:
$\preccurlyeq_{Q(\mathfrak{Q})}\, =\, \preccurlyeq_{\mathfrak{Q}}.$
Now the bijection $b_{St}$ is defined by
\[b_{St}\left(  \mathfrak{Q}\right)  =\left(  q_{_{P,Q\left(  \mathfrak{Q}
\right)  }},\Pi\left(  \mathfrak{Q}\, -\, q_{_{P,Q\left(  \mathfrak{Q}\right)  }}\right)  \right)  ,\]
where the reverse plane partition $\mathfrak{Q}\, -\, q_{_{P,Q\left(  \mathfrak{Q}\right)  }}$ is given by
\[\left(  \mathfrak{Q}\, -\, q_{_{P,Q\left(  \mathfrak{Q}\right)  }}\right)  \left(
\alpha\right)  =\mathfrak{Q\left(  \alpha\right)}\, -\, q_{_{P,Q\left(
\mathfrak{Q}\right)  }}\left(  \alpha\right)  ,\ \alpha\in\lambda.\]

To construct the inverse bijection, $b_{St}^{-1},$ we first associate to every standard Young tableau $Q\in\mathfrak{st}_{\lambda}$ and every partition
$p=(\Lambda_{1},\Lambda_{2},\dots,\Lambda_{n})$ the reverse plane partition $\mathfrak{Q}_{Q,p}$ of shape $\lambda,$ by
\[\mathfrak{Q}_{Q,p}\left(  i,j\right)  =\Lambda_{Q\left(  i,j\right)  }.\]
Then
\[b_{St}^{-1}\left(  q_{_{P,Q}},p\right)  =q_{_{P,Q}}+\,\mathfrak{Q}_{Q,p}.\]
The proof is finished. \hfill$\square$

The identity $\left(  \ref{01}\right)  $ is of independent interest. The principal specialization, $x_{i}\mapsto x^{i}$, turns
the `wrong' Schur function $\bar{s}_{\lambda}\left(  \mathbf{x}\right)$ into the generating function $\sigma_{\lambda}(x)$ for the number of reverse
plane partitions of shape $\lambda$, given by the Mac-Mahon--Stanley formula
\[\sigma_{\lambda}(x)=\frac{1}{\prod_{\alpha\in\lambda}\left(  1-x^{h_{\alpha} }\right)  },\]
where $h_{\alpha}$ is the hook length of a node $\alpha$ of $\lambda.$ The term $\bar{s}_{(n)}\left(  \mathbf{x}\right)$
turns into the generating function of Young diagrams with at most $n$ rows:
\[\sigma_{(n)}(x)=\frac{1}{\left(  1-x\right)  \left(  1-x^{2}\right)  ...\left(
1-x^{n}\right)  }.\]
Finally, the function $\mathfrak{h}_{\lambda}(\mathbf{x})  $ turns into the generating polynomial
$\pi_{\lambda}\left(  x\right)  $ of the sequence $\{ p_{\lambda,k}\}$ with $p_{\lambda,k}$ the number of $P$-pedestals
of volume $k\geq0$ for some $P \in\mathfrak{st}_{\lambda}$. We obtain
\begin{equation}\frac{1}{\prod_{\alpha\in\lambda}\left(  1-x^{h_{\alpha}}\right)  }=\frac{\pi_{\lambda}\left(  x\right)  }{\left(  1-x\right)  \left(
1-x^{2}\right)  ...\left(  1-x^{n}\right)  }. \label{04}\end{equation}
It follows from $\left(  \ref{02}\right)$ that the function
\begin{equation}\pi_{\lambda}\left(  x\right)  =\sum_{Q\in\mathfrak{st}_{\lambda}
}x^{ \vert\, q_{_{P,Q}}\vert} \label{03}\end{equation}
does not depend on the choice of $P$ while the contribution of an individual standard Young tableau $Q$ does.

The formula $(\ref{04})$ can be found in \cite{St}, but there the polynomial $\pi_{\lambda}\left(  x\right)  $ is given by any of two other
expressions:
\[\pi_{\lambda}\left(  x\right)  =x^{-l\left(  \lambda\right)  }\sum_{Q\in\mathfrak{st}_{\lambda}}x^{\mathrm{maj}\left(  Q\right)  }
\text{ and}\pi_{\lambda}\left(  x\right)  =x^{-l \left(  \lambda\right)  }\sum_{Q\in\mathfrak{st}_{\lambda}}x^{\mathrm{comaj}\left(  Q\right)  },\]
where $l\left(  \lambda\right) =\sum_{(i,j)\in\lambda}(i-1) $. It is interesting to note that in general neither of the two functions on
$\mathfrak{st}_{\lambda}$, $\mathrm{maj}(\cdot)  -l (\lambda)  $ and $\mathrm{comaj}(\cdot)  -l (\lambda)  $,  nor their partners for the transposed
to $\lambda$ Young diagram, belong to our family $\left\{  \mathrm{vol}\left(  q_{_{P,\ast}}\right)  :P\in\mathfrak{st}_{\lambda}\right\}$.
For example, take $\lambda=(3,2,1)$.

\end{document}